\documentstyle[twoside]{article}

\begin{document}
\
\begin{center}

{\LARGE{\bf{Duality and Boundary Yangians}}\footnote{ A talk presented 
at the International Workshop 
"Symmetries in Qantum Mechanics and Quantum Optics", Burgos 2-24 
September 1998. 
}\footnote{Supported in part by the Russian Foundation for
Fundamental Research grant N 97-01-01152 and by
the  Ministerio de Educacion y Cultura d'Espania, grant N SAB1995-0610.
}}

\vskip1cm

{\sc Lyakhovsky V.D.$^{3,4}$ and Ananikian D.N.$^4$  }

\vskip0.5cm

{\it $^3$ Departamento de Fisica Theorica, Facultad de Sciencias,
Universidad de Valladolid, E-47011, Valladolid, Spain 
}

{\it $^4$ Department of Theoretical Physics, St. Petersburg State
University, 198904, St. Petersburg, Russia}
\vskip0.3cm

\end{center}
\bigskip

\begin{abstract}
The existence of dual structures in a Yangian $Y(g)$ signify
that the latter belongs to multidimensional naturally 
parametrized variety of Hopf algebras. These varieties have
boundaries containing Yangians $Y(a)$ inequivalent to the     
original $Y(g)$. The new basic algebra $a$ is an algebra of the
cotangent bundle attributed to the dual structure. We show how 
to construct such boundary Yangians and study some of their properties.
We prove that the Hopf algebra $Y(a)$ is a quantization of 
a parametric solution of the classical Yang-Baxter equation.
The limiting procedure and the properties of the boundary
Yangians are demonstrated explicitly for the case of $Y(sl(2))$.

\end{abstract}

\section{Introduction}
The solutions of Yang-Baxter equation (YBE) depending on spectral parameter
are of special importance in mathematical physics. They are connected
with Wess-Zumino-Witten models and affine Toda field theories
\cite{BERN},\cite{BERN1}. The algebraic base of these theories is formed 
by the universal enveloping algebra $U(g[u])$ of polynomial current 
algebra $g[u]$ with $g$ being a finite-dimensional complex Lie algebra.

Yangians $Y(g)$ were introduced by Drinfeld \cite{DRIN1}  as a quantum 
deformations of the algebra $U(g[u])$.  
They correspond to the rational   
solutions of the classical Yang-Baxter equation (CYBE) first found by 
Belavin \cite{BEL} and 
completely described by Belavin and Drinfeld \cite{BELDR}.
The classification of rational solutions of CYBE for simple Lie algebras 
was performed by Stolin \cite{STO}.

The quantum deformation ${\cal R}(\lambda)$ of such rational solution 
\begin{equation}
\label{rat-sol}
r(\mu - \nu) = \frac{C^2}{\mu - \nu}
\end{equation}
(with $C^2$ being the second order Casimir element of the algebra g) 
satisfies the parametric Yang-Baxter equation,
$$
{\cal R}_{12}(\lambda_1 - \lambda_2){\cal R}_{13}(\lambda_1 - \lambda_3)
{\cal R}_{23}(\lambda_2 - \lambda_3) =
{\cal R}_{23}(\lambda_2 - \lambda_3){\cal R}_{13}(\lambda_1 - \lambda_3)
{\cal R}_{12}(\lambda_1 - \lambda_2),
$$
and realizes the morphism to the Yangian with the opposite 
comultiplication:
\begin{equation}
\label{pqybe}
\left( T_{\lambda} \otimes {\rm id} \right) \Delta^{\rm op}(a) =
{\cal R}(\lambda) \left( (T_{\lambda} \otimes {\rm id})
\Delta(a) \right) \left( {\cal R}(\lambda) \right)^{-1},
\end{equation}
where $T_{\lambda}$	is a parameter shifting operator. 

In this report we consider the problem of constructing
the quantizations of  $U(a[u])$ where the algebra $a$ is not semisimple
and also present rational solutions of CYBE for such algebras.
        
Our approach is based on studying the properties of the boundaries of
parametrized sets of Yangians. We find that under certain conditions the
corresponding algebraic constructions survive on the boundaries of   
the parametrized domain and can be explicitely described. To reach the
boundaries we use the restricted limiting procedure based on the
existence of dual structure in the corresponding Hopf algebra $H$.
The existence of a dual structure in $H$ means that there exists the
two-dimensional set of Hopf algebras containing $H$. One of the
properties of such sets is that when the corresponding parameters go to
its limiting values the algebraic construction survives. In most of the
cases the Hopf algebra thus obtained is inequivalent to the
original one. This was
demonstrated for deformation quantizations of the finite-dimensional
Lie algebras \cite{LYATKA}. Now we show that the same is true also for 
Yangians and other possible  quantizations of  $U(g[u])$.

\section{Yangian $Y(sl(2))$ and its natural limits}

The Yangian $Y(sl(2))$ is a Hopf algebra generated by the elements
$ \{ e_k, h_k, f_k \} $ $ (k \in {\bf Z}_+ )  $ with relations
\begin{equation}
\label{ini-my}
\begin{array}{l}  
\begin{array}{lcl}
\, [h_k , h_l ]  =  0, & & [e_k , f_l ]  =  h_{k+l},\\[1mm]
\, [h_0 , e_l ]  =  2e_{l}, & & [h_0 , f_l ]  =  -2f_{l},\\[1mm]
\end{array} \\
\begin{array}{lcl}
\, [h_{k+1} , e_l ] -[h_k , e_{l+1} ] & = &  \hbar \{ h_k,e_l \},\\[1mm]
\, [h_{k+1} , f_l ] -[h_k , f_{l+1} ] & = &  -\hbar \{ h_k,f_l \},\\[1mm]
\, [e_{k+1} , e_l ] -[e_k , e_{l+1} ] & = &  \hbar \{ e_k,e_l \},\\[1mm] 
\, [f_{k+1} , f_l ] -[f_k , f_{l+1} ] & = &  -\hbar \{ f_k,f_l \},\\
\end{array}
\end{array}
\end{equation}
where $\hbar$ is the deformation parameter. The coproducts for the
generators of $sl(2) \in Y(sl(2))$ rest primitive:
\begin{equation}
\label{ini-py}  
\Delta (x) = x \otimes 1 + 1 \otimes x; \, \, \, x \in sl(2)
, 
\end{equation}
while the nontrivial coalgebraic part is uniquely defined by the
following comultiplications (and the multiplications (\ref{ini-my})
 above):
\begin{equation}
\label{ini-cy}  
\begin{array}{lcl}
\Delta (e_1) & = &  e_1 \otimes 1 + 1 \otimes e_1 + \hbar h_0 \otimes
 e_0, \\
\Delta (f_1) & = &  f_1 \otimes 1 + 1 \otimes f_1 + \hbar f_0 \otimes
 h_0. \\
\end{array}
\end{equation}
To be able to study general properties of Yangian's Hopf structure
the generating functions formalism is especially convenient. In terms of

\begin{equation}
\label{gen-fs}
\begin{array}{l}
e(u) := \sum_{k \geq 0} e_k u^{-k-1}, \;\;
f(u) := \sum_{k \geq 0} f_k u^{-k-1}, \\[1mm]
h(u) :=1 + \hbar\chi(u) := 1 + \hbar \sum_{k \geq 0} h_k u^{-k-1}
\end{array}
\end{equation} 
the compositions generated by (\ref{ini-my},\ref{ini-py},\ref{ini-cy})
look like
\begin{equation}  
\label{ini-myg}   
\begin{array}{lcl}
\,[h(u) , h(v)]  =  0, & & [e(u), f(v)]  =    
- \frac{1}{\hbar}\frac{h(u) - h(v)}{u-v} ,\\[2mm]
\,[h(u), e(v)]  =                         
- \hbar\frac{ \{ h(u),(e(u) - e(v))\} }{u-v} ,
 & & [h(u), f(v)]  = 
  \hbar\frac{ \{ h(u),(f(u) - f(v))\} }{u-v} ,\\[2mm]
\,[e(u), e(v)]  =                       
- \hbar\frac{ (e(u) - e(v))^2 }{u-v} ,
& & 
[f(u), f(v)]  =                      
 \hbar\frac{ (f(u) - f(v))^2 }{u-v} .
\end{array}   
\end{equation}
The coproducts for the generating functions are written in the
form proposed by Molev \cite{MOL}
\begin{equation}  
\label{ini-cyg}   
\begin{array}{lcl}
\Delta (e(u)) & = & e(u) \otimes 1 + \sum^{\infty}_{k=0} 
(-1)^k \hbar^{2k} (f(u + \hbar))^k h(u) \otimes (e(u))^{k+1},\\[3mm]
\Delta (f(u)) & = & 1 \otimes f(u) + \sum^{\infty}_{k=0} 
(-1)^k \hbar^{2k} (f(u))^{k+1} \otimes h(u)(e(u + \hbar))^k,\\[3mm]
\Delta (h(u)) & = & \sum^{\infty}_{k=0}          
(-1)^k (k+1) \hbar^{2k} (f(u + \hbar))^{k} h(u) \otimes h(u)(e(u
+ \hbar))^{k}.
\end{array}   
\end{equation}

The dual structure in $Y(sl(2))$ is inherited from the basic  
subalgebra $sl(2)$ interpreted as a classical double. When the
double structure is in its prime form the dual parameters and
the corresponding analytical family of Hopf algebras can be
canonically introduced (see \cite{LYM}). 
In the case of Yangian this double is factorized. To get the   
necessary parametrization we must reconstruct the prefactorized
relations. Such reconstruction can be achieved if one takes
into account that the unfactorized classical double of the
initial $sl(2)$ subalgebra has the form
\begin{equation}
\label{c-dub}
\begin{array}{lcl}
\,[h, h'] = 0, & & [e, f] = \frac{1}{2} (h + h'), \\
\,[h, e] = 2 e, & & [h, f] = -2 f, \\
\,[h', e] = 2 e, & & [h', f] = -2 f.
\end{array}   
\end{equation}
This means that the transformation
\begin{equation}
\label{subs}
\begin{array}{lcl}
e_k & \rightarrow & \frac{1}{p} e_k, \\[1mm]     
f_k & \rightarrow & \frac{1}{t} f_k, \\[1mm]     
h_k & \rightarrow & \frac{1}{2}(\frac{h_k}{p} + \frac{h'_k}{t})
\end{array}   
\end{equation}
accompanied by the rescaling of the original $sl(2)$ structure
constants will lead to the parametrization that would not be  
canonical in the whole two-dimensional domain but might work  
well just in the neighborhood of its boundaries. According to 
the general rules \cite{LYM} the deformation parameter must be
also rescaled:  
\begin{equation}
\hbar \rightarrow pt.
\end{equation}
All these transformations and rescalings produce the pa\-ra\-met\-rized
algebraic construction  $Y_{pt}(sl(2))$ well defined only when
one of the parameters ($p$ or $t$) is small.

When a deformation quantization algebra admits a canonical parametrization
(according to some dual structure) its quasiclassical limits can be
found. In our case the limits $ Y_{p,0}(a) := 
\lim_{t \rightarrow 0}Y_{pt}(sl(2))$ 
and $Y_{0,t}(b) := \lim_{p \rightarrow 0}Y_{pt}(sl(2))$ are 
equivalent due to the 
selfduality of the Borel subalgebra in the classical double
algebra (\ref{c-dub}). Thus we need to consider only one of
these algebraic constructions, that we call boundary Yangians. 
The limiting procedure
gives the following structure constants for the boundary Yangian
$Y_{p,0}(a)$:
\begin{equation}
\label{lim-my}
\begin{array}{l}
\begin{array}{lll}
\, [h_k , h_l ]  =  0, &  [h'_k , h'_l ]  =  0,
&  [h_k , h'_l ]  =  0,\\[1mm]
\, [h_0 , e_l ]  =  4pe_{l}, &  [h_0 , f_l ]  =  -4pf_{l}, & \\[1mm]
\, [h'_0 , e_l ]  = 0, &  [h'_0 , f_l ]  =  0,
& [e_k , f_l ]  =  \frac{p}{2} h'_{k+l}, \\[1mm]
\end{array}\\[2mm]
\begin{array}{lll}
\, [h_{k+1} , e_l ] -[h_k , e_{l+1} ] & = &  p^2 \{ h'_k,e_l \},\\[1mm] 
\, [h_{k+1} , f_l ] -[h_k , f_{l+1} ] & = &  -p^2 \{ h'_k,f_l \},\\[1mm]
\, [h'_{k+1} , e_l ] -[h'_k , e_{l+1} ] & = &  0,\\[1mm] 
\, [h'_{k+1} , f_l ] -[h'_k , f_{l+1} ] & = &  0,\\[1mm]
\, [e_{k+1} , e_l ] -[e_k , e_{l+1} ] & = &  0,\\[1mm]
\, [f_{k+1} , f_l ] -[f_k , f_{l+1} ] & = &  0,\\[1mm]   
\end{array}
\end{array}
\end{equation}
The coproducts for the "zero mode"
generators $e_0,f_0,h_0,h'_0$ rest primitive. The others are
defined by the relations:
\begin{equation}
\label{lim-cy}  
\begin{array}{lcl}
\Delta (e_1) & = &  e_1 \otimes 1 + 1 \otimes e_1 + \frac{p}{2}
h'_0 \otimes  e_0, \\[1mm]
\Delta (f_1) & = &  f_1 \otimes 1 + 1 \otimes f_1 + \frac{p}{2}
f_0 \otimes h'_0, \\[1mm]
\Delta (h_1) & = &  h_1 \otimes 1 + 1 \otimes h_1 + \frac{p}{2}
(h'_0 \otimes h_0 + h_0 \otimes h'_0) \\[1mm]
  & &  - 4p  f_0 \otimes e_0,\\[1mm]
\Delta (h'_1) & = &  h'_1 \otimes 1 + 1 \otimes h'_1
+ \frac{p}{2}(h'_0 \otimes h'_0). \\[1mm]
\end{array}
\end{equation}
The internal structure of this Hopf algebra becomes more
transparent in terms of generating functions:
\begin{equation}
\label{lim-myg} 
\begin{array}{l}
\begin{array}{lll}
\, [\chi(u) , \chi(v) ]  =  0, &  [\chi'(u) , \chi'(v) ]  =  0,
&  [\chi(u) , \chi'(v) ]  =  0,\\[1mm]
\, [\chi'(u) , e(v) ]  = 0, &  [\chi'(u) , f(v) ]  =  0, & \\[2mm]
\, [e(u) , f(v) ]  = 
- \frac{p}{2}\frac{ \chi'(u) - \chi'(v) }{u-v}, \\[1mm]
\end{array}\\[2mm]

\begin{array}{lll}
\, [\chi(u) , e(v) ] & = &  - \frac{p}{u-v}
\{ 2+p \chi'(u), e(u) - e(v)  \} \\[1mm]
\, [\chi(u) , f(v) ] & = &  \frac{p}{u-v}
\{ 2+p \chi'(u), f(u) - f(v)  \}  \\[1mm]
\end{array}\\[2mm]
\end{array}
\end{equation}

\begin{equation}
\label{lim-cyg}
\begin{array}{lcl}
\Delta (e(u)) & = &  e(u) \otimes 1 + 1 \otimes e(u) + \frac{p}{2}
\chi'(u) \otimes  e(u), \\[1mm]
\Delta (f(u)) & = &  f(u) \otimes 1 + 1 \otimes f(u) + \frac{p}{2}
f(u) \otimes \chi'(u), \\[1mm]
\Delta (\chi'(u)) & = &  \chi'(u) \otimes 1 + 1 \otimes \chi'(u)
+ \frac{p}{2}(\chi'(u) \otimes \chi'(u)), \\[1mm]
\Delta (\chi(u)) & = &  \chi(u) \otimes 1 + 1 \otimes \chi(u) + 
\\[1mm]
& & +\frac{p}{2}(\chi'(u) \otimes \chi(u) + 
   \chi(u) \otimes \chi'(u)) \\[1mm]
  & &  - 4p  f(u) (1 + \frac{p}{2} \chi'(u)) \otimes
 (1 + \frac{p}{2} \chi'(u)) e(u),\\[1mm]
\end{array}
\end{equation}
where $e(u), f(u)$ and $\chi(u)$ are as in (\ref{gen-fs}) while $\chi'(u)$ 
is the analog of $\chi(u)$ for the generator $h'$. 
This Hopf algebra $Y_{p,0}(a)$ is a Yangian for a nonsemisimple
Lie algebra $a$ with the following compositions:
\begin{equation}
\label{lim-g}
\begin{array}{c}
\, [h'_0, h_0] = 0,\\[1mm]
\begin{array}{ll}
\, [h_0, e_0] =4pe_0, & [h'_0, e_0] =0, \\[1mm]
\, [h_0, f_0] =-4pf_0, & [h'_0, f_0] =0, \\
\end{array}\\[2mm]
\, [e_0, f_0] =\frac{p}{2}h'_0,
\end{array}
\end{equation}
This is just the cotangent bundle algebra for the two-dimensional
Borel subalgebra of sl(2).

 For possible applications it is crutially important to know 
whether the Yangians like $Y_{p,o}(a)$ are pseudotriangular or not,
that is whether they have the universal ${\cal R}$-matrix providing the
property (\ref{pqybe}).
Below we shall present
the indications that such ${\cal R}$-matrises exist.

The direct way to solve the problem is to construct the
parametrized version of the expansion terms for 
 ${\cal R}$-matrix in the general case of $Y_{p,t}(sl(2))$ and
to check their limits. In our case this aproach doesn't work. The first
nontrivial term -- the classical $r$-matrix after being parametrized
according to the transformation (\ref{subs}) diverges in the limit
point. As we have seen above the limiting procedure must be 
accompanied by the rescaling of commutators. This clearly indicates
that certain terms of the $r$-matrix must be rescaled in the 
neighborhood of the limit point. It is not difficult to find these
terms and to perform an adequate rescaling. The result is formulated 
in the statement that follows.

\underline{Theorem.} {\sl The boundary Yangian $Y_{p,0}(a)$ (
correspondingly $Y_{0,t}(a))$ originating from $Y(sl(2))$ is the 
deformation 
quantization of the polinomial current algebra based on the 
algebra  $a$ (\ref{lim-g}). The first order expansion term for  
this deformation is defined by the following solution of the 
classical Yang-Baxter equation for the algebra $a$:
\begin{equation}
\label{rat-sol2}
r(u,v) = \frac{1}{u-v} \left[ \frac{1}{8}(h_0 \otimes h'_0 + 
h'_0 \otimes h_0) + e_0 \otimes f_0 + f_0 \otimes e_0 
		       \right].
\end{equation}
}        

The first assertion becomes obvious when the classical limit of
$Y_{p,0}(a)$ is considered. The validity of the last statement
can be checked by the direct computation of the dual Lie algebra 
defined by (\ref{rat-sol2}); it coinsides with the co-Lie structure
that one observes in (\ref{lim-cy}) or (\ref{lim-cyg}).

To make the demonstration most transparent we have studied the 
simpliest possible case -- the Yangian based on $sl(2)$ algebra.
As a result our boundary Yangian $Y(a)$ admits some additional 
simplifications. The subalgebra generated by $\chi'(u)$ forms a Hopf
ideal  $J(\chi'(u)) \in Y(a)$. In the factor algebra 
$\frac{Y(a)}{J(\chi'(u))} \equiv Y(c) $
the only nontrivial relations are 
\begin{equation}
\label{fact-m}
\begin{array}{lll}
\, [\chi(u) , e(v) ] & = &  - \frac{4p}{u-v}
(e(u) - e(v)), \\[1mm]
\, [\chi(u) , f(v) ] & = &  \frac{4p}{u-v}
(f(u) - f(v))  
\end{array}
\end{equation}
and
\begin{equation}
\label{fact-c}  
\Delta (\chi(u))  =   \chi(u) \otimes 1 + 1 \otimes \chi(u) 
   - 4p f(u) \otimes e(u).
\end{equation}
This Yangian $Y(c) \approx U_{\cal F}(c[u])$ is a 
quantized algebra of currents 
for the Lie algebra $c$, where $c$ is 
the algebra $(sl(2))^{\rm contr}$ with 
the trivially contracted composition $ [e,f]=0 $.

The defining relations show that the multiplications in $Y(c)$ are 
undeformed; the commutators are classical and coinside with those of
$U(c [ u ] )$,
\begin{equation}
\label{fact-m-g}
\begin{array}{llllll}
\, [h_k,h_l] & = & 0; & [h_0,e_l] & = & 4pe_l ; \\
\, [e_k,f_l] & = & 0; & [h_0,f_l] & = & -4pf_l ; \\
\, [e_{k+1},e_l] & = & [e_{k},e_{l+1}]; 
& [f_{k+1},f_l] & = & [f_{k},f_{l+1}] ;\\
\, [h_{k+1},e_l] & = & [h_{k},e_{l+1}]; 
& [h_{k+1},f_l] & = & [h_{k},f_{l+1}] .
\end{array}
\end{equation}  
It also has trivial coproducts 
for $e_k$, $f_k$ \, $(k \in {\bf Z}^+)$ and $h_0$.
The nontrivial costructure in $Y(c)$ is generated by the relation
$$
\Delta (h_1) = h_1 \otimes 1 + 1 \otimes h_1 - 4pf_0 \otimes e_0 
$$
and the compositions (\ref{fact-m-g}).

This quantized algebra can be obtained from $U(c[u])$ by a twisting
procedure: $U(c[u]) \stackrel{\cal F}{\rightarrow} Y(c)$. The carrier
of this twist  is an abelian subalgebra generated by the primitive 
elements $e_0$ and $f_0$ and the twisting element has the form
$$
{\cal F}(u-v) = \exp{\left(  \frac{1}{v - u} f_0 \otimes e_0 \right)}.
$$
This gives the possibility to write down the universal element for $Y(c)$,
$$
{\cal R}(u-v) = \exp{ \left( \frac{1}{u-v} ( f_0 \otimes e_0 
+ e_0 \otimes f_0 ) \right) }.  
$$
It can be checked that this ${\cal R}$-matrix 
provides the pseudotriangularity
of $Y(c)$. The first nontrivial term of the power series expansion of
${\cal R}$ coinsides with 
the classical $r$-matrix that can be obtained from the expression
(\ref{rat-sol2}) after the factorization by the ideal 
$J(h'_k)|_{k \in {\bf Z}^{+}}$.

\section{Conclusions}

We have demonstrated that dual structures signify the existence 
of the nontrivial limiting algebraic objects corresponding to the
rational solutions of CYBE based on nonsemisimple Lie algebras. 
Thus rational solutions of the type (\ref{rat-sol}) 
are not necessarily connected with the nondegeneracy of the Killing form.

The existence of such rational solutions for the so called symmetric 
algebras was established in \cite{STOBF}. It was also mentioned there that 
probably Yangians for symmetric algebras exist. 
Our result proves this supposition. An important subclass of 
symmetric algebras is presented by Manin triples. These are the algebras
$(a, b(2), b(2)^*)$ (with two-dimensional Borel algebra $b(2)$)  
that form a Manin triple in our case and its  
 characteristic nondegenerate form was used to construct the 
invariant element in (\ref{rat-sol2}).

The procedure proposed above is quite general and can be applied 
to any Hopf algebra that have the form of deformation quantization.
For example, it can be applied to the Yangian double (DY) where the
dual structure is more simple than in the case of the Yangian itself
-- the double is not factorized in DY. It must be noted that in this case 
the situation with the canonical limits for the  ${\cal R}$-matrix 
is complicated by the specific form of dualization established for DY
\cite{KHOR}. It depends on the   
quantization parameter and fails in the canonical limits.
To reobtain quasitriangularity for the corresponding boundary Yangian one 
must  reformulate the structures of the
Yangian double in terms of the canonical dualization.

It is significant that after the factorization the boundary Yangian 
can be presented in the form of a twist so that its ${\cal R}$-matrix can
be written explicitly.
    
{\section*{Acknowledgments}}
\noindent
One of the authors (L.V.D.) expresses his sincere gratitude to the 
organizers of the International Workshop SQMQO in Burgos for their 
hospitality.

\end{document}